# Verification of Collatz Conjecture: An algorithmic approach based on binary representation of integers

Venkatesulu Mandadi[1] and Devi Paramwswari[2]

**Abstract**

Lothar Collatz had proposed in 1937 a conjecture in number theory called Collatz conjecture. Till today there is no evidence of proving or disproving the conjecture. In this paper, we propose an algorithmic approach for verification of the Collatz conjecture based on bit representation of integers. The scheme neither encounters any cycles in the so called Collatz sequence and nor the sequence grows indefinitely. Experimental results show that the Collatz sequence starting at the given integer , oscillates for finite number of times, never exceeds 1.7 times (scaling factor) size of the starting integer and finally reaches the value 1. The experimental results show strong evidence that conjecture is correct and paves a way for theoretical proof.

**Introduction:**

In 1937 Lothar Collatz proposed a conjecture which states that given any positive integer x, the function f(x) defined as x/2 if x is even and 3x+1 if x is odd generates a finite sequence called the Collatz sequence which eventually reaches the value. Then after it loops between the values 1 → 4 →2 →1, called a trivial cycle. In other words, the conjectures say that the Collatz sequence does not contain nontrivial cycles and also does not grow indefinitely. In the past and present many authors [1-11] try to establish the validity of the conjecture in affirmatively but till today there is no theoretical proof for the same and there is no counter example to disprove the conjecture as well.

In this paper, we propose an algorithmic approach based on binary representation of integers for verification of the Collatz conjecture. The scheme encounters no cycles in the so called Collatz sequence and observes that the Collatz sequence does not grows indefinitely. An experiment results show that the Collatz sequence starting at any given positive integer , oscillates for finite number of times, no element in the sequence ever exceeds 1.7 times (the scaling factor) size of the starting integer and finally reaches the value 1. In other words, for any given positive integer, the elements of Collatz sequence could grow to a size less than 1.7 times the size of the starting number and in finite number of steps, called the stopping time, the sequences eventually reaches the value 1. These results show strong evidence that conjecture is correct and paves a way for a possibility of theoretical proof.

**Verification Scheme:** The given positive integer is converted into binary form. Usually, the left most bit represents the most significant bit and the right most bit is the least significant bit. But for our scheme the binary representation is reversed and taken as the input for the verification algorithm. The binary representation of 2x is just ONE place right shift (note that we have reversed the form) of the binary form of x. The expression 3x + 1 are written x + 2x +1, all taken in binary form with the usual binary addition operations. For an even integer x, x/2 is just ONE place left shift of the binary expression of x.

**Algorithm**:

Step1: Without loss of generality, we may assume that x is an odd integer otherwise keep dividing the integer successively by 2 until we get an odd integer.

Step2: Convert the given integer x into binary form and reverse the binary representation.

Step3: Represent the expression 3x + 1 as binary multiplication and addition.

Step4: Do the division of an even integer by 2 by shifting the resulting binary representation by one place to the right.

Step5: Repeat steps3 and 4, appropriately, until step4 results in the value 1.

**Experimental Results:**

In this section we present experimental results on verification of Collatz conjecture. The algorithm is implemented using C Programming in a DELL system with Intel®Pentium®3558U @1.70 GHz processor with 4GB RAM. The graphs are drawn using R.

**Convert Binary Form:**

We have taken RSA-100: 1522605027922533360535618378132637429718068114961380688657908494580122963258952897654000350692006139 and converted into binary form. The decimal form contains 100 Decimal digits and its converted binary form has 330 .The binary list is reversed (as per our design, leftmost bit should be the least significant bit) and given as an input for step3 of the algorithm. The experiment result gave 1566 right shifts (even number & division by 2) and 780 loop iterations (odd integer & 3x+1) before the Collatz sequence

Below, we present experimental results for integers in binary form with all bits equal to 1. These represent the largest possible integer in the binary form with that size. Numbers of zeros represent division by 2 and loop iteration represent occurrence of an odd integer before the Collatz sequence reaches value 1.

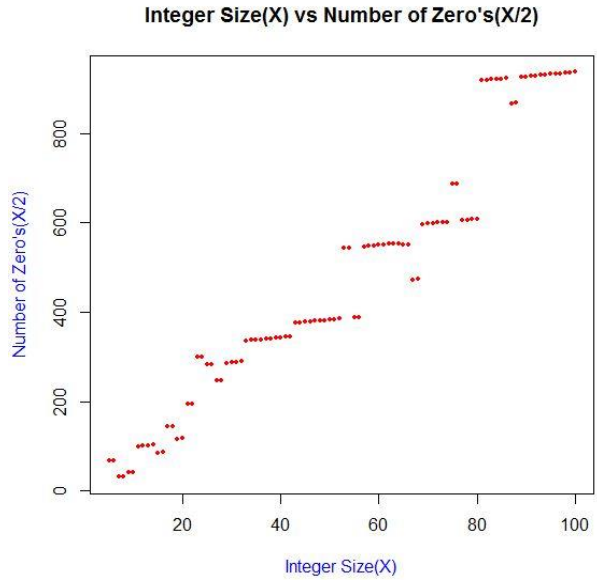

Fig 1: Integer Size(x) vs. Number of Zeros (x/2)

In fig.1, integer size(x) indicates the size of the numbers in binary form and x/2 indicates the number of times division by 2 is carried out before the Collatz sequence reaches the value 1

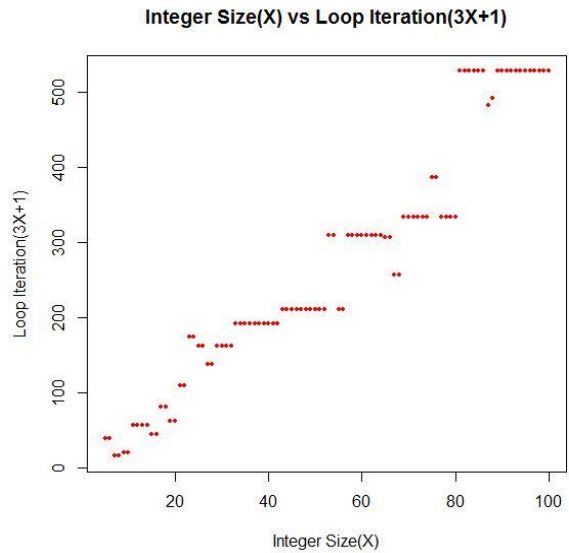

Fig 2: Integer Size(x) vs. Loop Iteration (3x+1)

In fig 2, integer size(x) denotes the size of the numbers in binary form and 3x+1 denotes the number of times the integer becomes odd during the process before the Collatz sequence reaches value1.

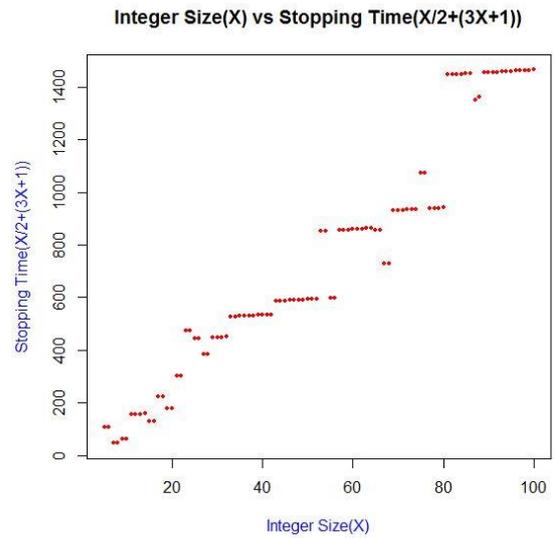

Fig 3: Integer Size(x) vs. Stopping time

In fig 3, integer size(x) denotes the size of the numbers in binary form and stopping time denotes the length of the Collatz sequence before it reaches value 1.

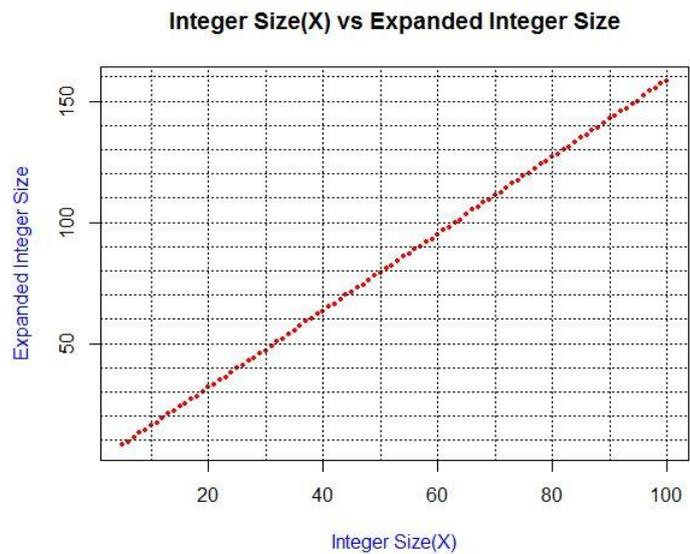

Fig 4: Integer Size(x) vs. Expanded Integer Size

In fig 4, integer size(x) denotes the size in binary form and expended integer size denotes the maximum size in binary form of any element of the Collatz sequence before it reaches value 1.

```
R Console
                85         86         87         88         89         90         91
140.959577 142.544086 144.128595 145.713104 147.297613 148.882122 150.466630
                92         93         94         95         96
152.051139 153.635648 155.220157 156.804666 158.389175
> summary(lm(V2~V1))

Call:
lm(formula = V2 ~ V1)

Residuals:
     Min       1Q   Median       3Q      Max
-0.48766 -0.24471 -0.00176  0.24118  0.49824

Coefficients:
             Estimate Std. Error  t value Pr(>|t|)
(Intercept) -0.061720   0.063821   -0.967    0.336
V1           1.584509   0.001075 1473.859   <2e-16 ***
---
Signif. codes:  0 '***' 0.001 '**' 0.01 '*' 0.05 '.' 0.1 ' ' 1

Residual standard error: 0.2919 on 94 degrees of freedom
Multiple R-squared:      1,     Adjusted R-squared:      1
F-statistic: 2.172e+06 on 1 and 94 DF,  p-value: < 2.2e-16

>
```

Fig 5: Regression Analysis

Here , we have considered integers with bits size from 5 bits to 100 bits and the maximum size of any element of any element in the Collatz sequence before it reaches value 1.The graph is drawn by applying Linear regression model in R. We got R-squared value is 1 which means that model is totally perfect.

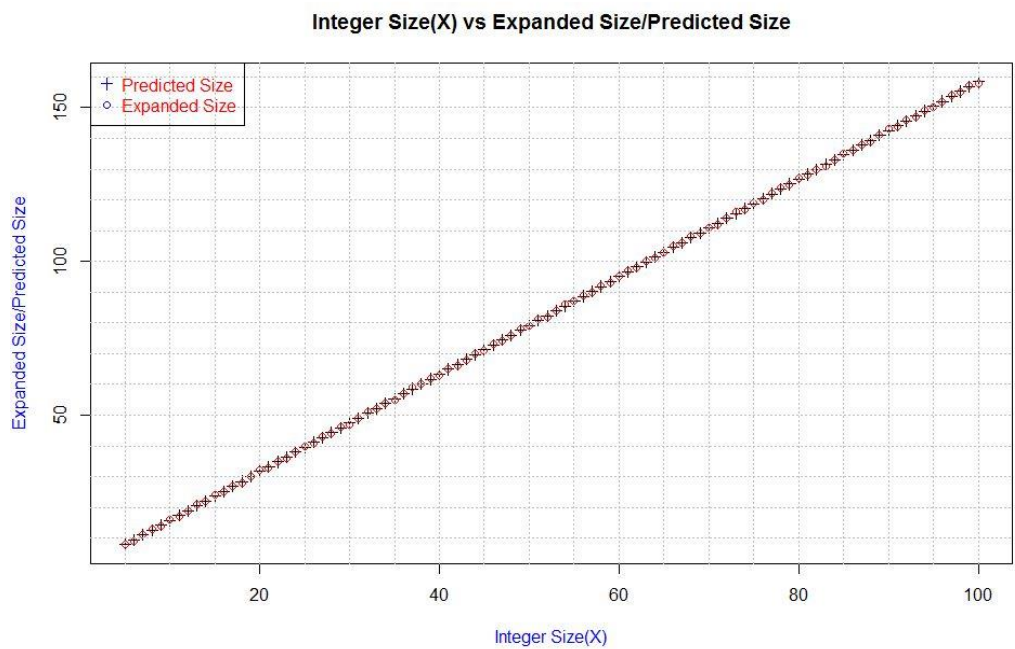

Fig 6: Integer Size vs. Expanded/Predicted Size

Fig.6 shows that Expended size and Predicted size of any element of the Collatz sequence are approximately the same.

We also performed the experiments for integers with binary representation from represented 100 to 3000 bits size.

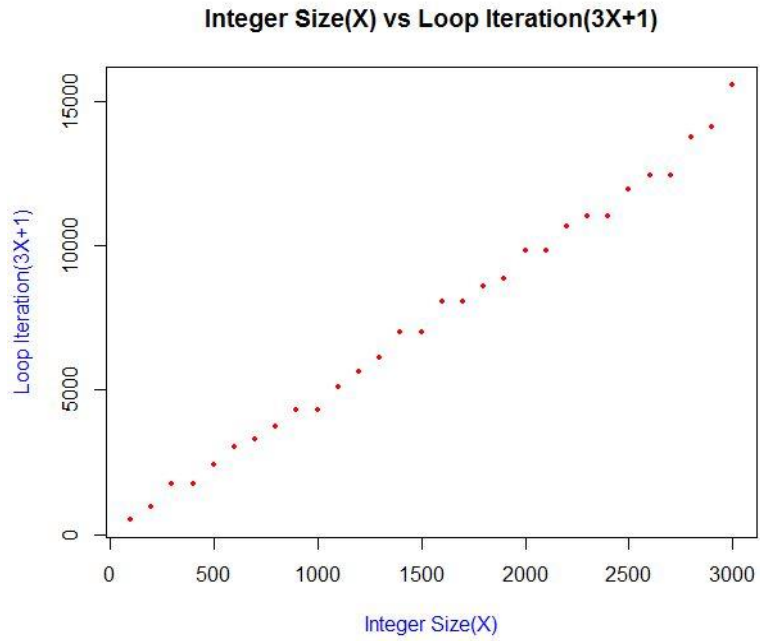

Fig 7: Integer Size(x) vs. Loop Iteration (3x+1)

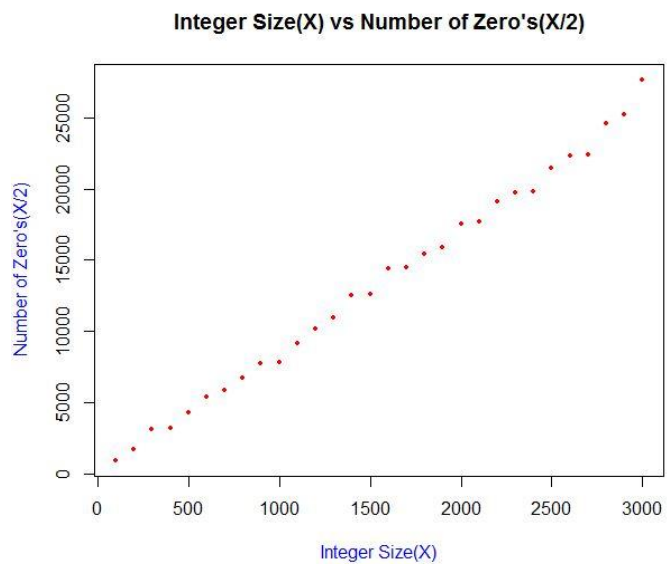

Fig 8: Integer Size(x) vs. Number of Zero's(x/2)

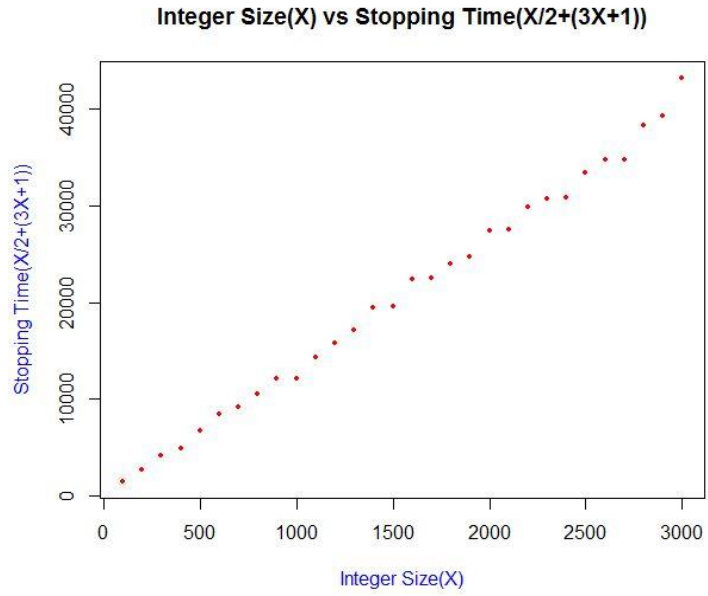

Fig 9: Integer Size(x) vs. Stopping Time(x/2) + (3x+1))

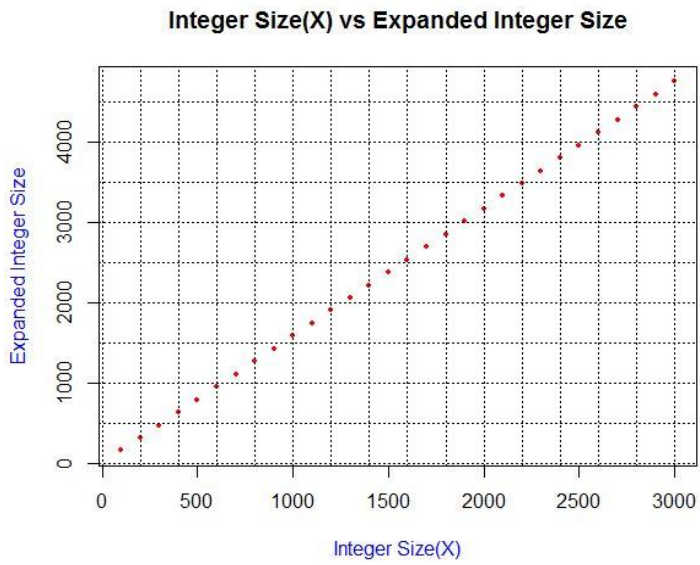

Fig 10: Integer Size(x) vs. Expanded Integer Size

```
R Console
      17         18         19         20         21         22         23         24
2694.2525 2852.7542 3011.2558 3169.7575 3328.2592 3486.7608 3645.2625 3803.7642
      25         26         27         28         29         30
3962.2659 4120.7675 4279.2692 4437.7709 4596.2725 4754.7742
> summary(lm(V2~V1))

Call:
lm(formula = V2 ~ V1)

Residuals:
    Min      1Q  Median      3Q     Max
-0.2725 -0.2483  0.0000  0.2483  0.2725

Coefficients:
              Estimate Std. Error  t value Pr(>|t|)
(Intercept) -2.759e-01  9.674e-02   -2.852  0.00808 **
V1           1.585e+00  5.449e-05 29086.262  < 2e-16 ***
---
Signif. codes:  0 '***' 0.001 '**' 0.01 '*' 0.05 '.' 0.1 ' ' 1

Residual standard error: 0.2583 on 28 degrees of freedom
Multiple R-squared:      1,    Adjusted R-squared:      1
F-statistic: 8.46e+08 on 1 and 28 DF,  p-value: < 2.2e-16

>
```

Fig 11: Regression Analysis

Here , we have considered integers with bits size from 100bits to 3000 bits and the maximum size of any element of any element in the Collatz sequence before it reaches value 1.The graph is drawn by applying Linear regression model in R. We got R-squared value is 1 which means that model is totally perfect. It is statistical measures for goodness-of-fit.

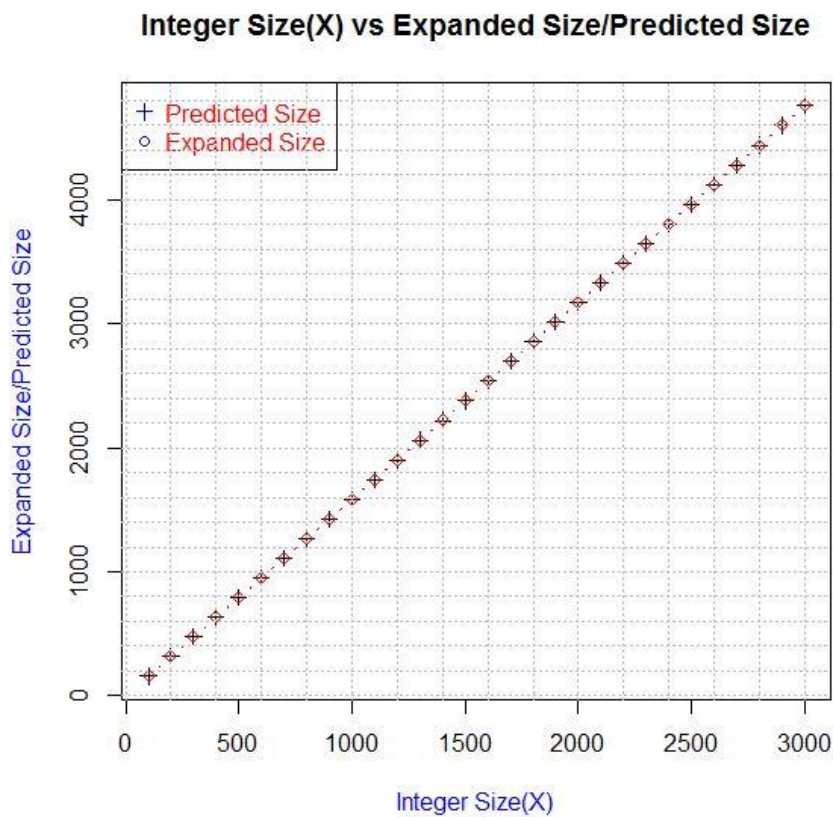

Fig 12: Integer Size vs. Expanded/Predicted Size

The table below shows the experimental results on verification of Collatz Conjecture where in we have introduced some zeros in the binary form. The second column denotes number zeros in the binary form of the given integer size, x/2 and 3x+1. Fourth column denotes the number times division by 2 happens and fifth column denotes the number times the odd integer is encountered during the Collatz sequence generation before it reaches value 1. The third column denotes the maximum size of any element of the Collatz sequence.

| INTEGER SIZE | NUMBER OF ZERO'S | EXPANDED INTEGER SIZE | x/2 | 3x+1 |
|---|---|---|---|---|
| 100 | 1 | 155 | 937 | 528 |
| 100 | 2 | 116 | 588 | 309 |
| 100 | 3 | 112 | 429 | 210 |
| 100 | 4 | 112 | 708 | 386 |
| 100 | 5 | 114 | 431 | 210 |
| 200 | 1 | 245 | 1102 | 569 |
| 200 | 2 | 244 | 1316 | 704 |
| 200 | 3 | 244 | 1140 | 593 |
| 200 | 4 | 244 | 1316 | 704 |
| 200 | 5 | 244 | 1102 | 569 |
| 300 | 1 | 312 | 1443 | 721 |
| 300 | 2 | 312 | 1454 | 728 |
| 300 | 3 | 312 | 2020 | 1085 |
| 300 | 4 | 312 | 1454 | 728 |
| 300 | 5 | 315 | 2004 | 1075 |
| 400 | 1 | 547 | 3212 | 1774 |
| 400 | 2 | 547 | 3212 | 1774 |
| 400 | 3 | 545 | 2811 | 1521 |
| 400 | 4 | 545 | 3212 | 1774 |
| 400 | 5 | 545 | 2765 | 1492 |
| 500 | 1 | 672 | 3648 | 1986 |
| 500 | 2 | 676 | 4030 | 2227 |
| 500 | 3 | 676 | 3277 | 1752 |
| 500 | 4 | 676 | 3648 | 1986 |
| 500 | 5 | 677 | 3277 | 1752 |

| 1000 | 1 | 1246 | 6370 | 3388 |
| --- | --- | --- | --- | --- |
| 1000 | 2 | 1247 | 6266 | 7892 |
| 1000 | 3 | 1246 | 6434 | 7998 |
| 1000 | 4 | 1246 | 6266 | 7892 |
| 1000 | 5 | 1247 | 6266 | 7892 |
| 1500 | 1 | 2051 | 10982 | 10552 |
| 1500 | 2 | 2049 | 10982 | 10552 |
| 1500 | 3 | 2047 | 11399 | 10815 |
| 1500 | 4 | 2047 | 10982 | 10552 |
| 1500 | 5 | 2046 | 10981 | 10552 |
| 2000 | 1 | 2639 | 14839 | 12670 |
| 2000 | 2 | 2639 | 14270 | 12311 |
| 2000 | 3 | 2639 | 13709 | 11957 |
| 2000 | 4 | 2639 | 13159 | 11610 |
| 2000 | 5 | 2639 | 13159 | 11610 |
| 2500 | 1 | 3506 | 18192 | 14470 |
| 2500 | 2 | 3505 | 18417 | 14612 |
| 2500 | 3 | 3504 | 18417 | 14612 |
| 2500 | 4 | 3504 | 18417 | 14612 |
| 2500 | 5 | 3504 | 18875 | 14901 |
| 3000 | 1 | 4379 | 23239 | 17339 |
| 3000 | 2 | 4378 | 24298 | 18007 |
| 3000 | 3 | 4378 | 23672 | 17612 |
| 3000 | 4 | 4378 | 23198 | 17313 |
| 3000 | 5 | 4379 | 23198 | 17313 |

**Conclusion:** We have presented a new algorithmic approach for Collatz conjecture verification based on binary representation, multiplication, addition and division by 2, all done in binary set up. We observe that there are NO cycles in Collatz sequence and no element of the sequence exceeds in size 1.7 times the size of the given (starting) number. We have verified the results with all entries in binary form equal to 1, being the largest possible integer with that size. Due to limit in available computing resources, we verified the conjecture up to binary string of size 3000. Given, enough computing resources, our scheme can verify the conjecture for any given integer however big it may be.

**Acknowledgement:** We thank Professor Dr.Sampath, Department of Computer Science and Information Technology, Kalasalingam Academy of Research and Education and Mr. Krishna Kumar, Bangalore for their valuable inputs in representing the experimental results through graphs. The second author thanks the management of Kalasalingam Academy of Research and Education for providing scholarship to carry out the research

[1]Department of Information Technology
[2]Department of Computer Science and Information Technology
Kalasalingam Academy of Research and Education
Anandnagar, Krishnankoil-626126, Srivilliputtur (Taluk),
Virudhunagar (Dist.), Tamilnadu, India
Email: m.venkatesulu@klu.ac.in
https://orcid.org/0000-0002-3972-9929